# ON MODELING THE RESPONSE OF SYNOVIAL FLUID: UNSTEADY FLOW OF A SHEAR-THINNING, CHEMICALLY-REACTING FLUID MIXTURE

CRAIG BRIDGES, SATISH KARRA, AND K. R. RAJAGOPAL

ABSTRACT. We study the flow of a shear-thinning, chemically-reacting fluid that could be used to model the flow of the synovial fluid. The actual geometry where the flow of the synovial fluid takes place is very complicated, and therefore the governing equations are not amenable to simple mathematical analysis. In order to understand the response of the model, we choose to study the flow in a simple geometry. While the flow domain is not a geometry relevant to the flow of the synovial fluid in the human body it yet provides a flow which can be used to assess the efficacy of different models that have been proposed to describe synovial fluids. We study the flow in the annular region between two cylinders, one of which is undergoing unsteady oscillations about their common axis, in order to understand the quintessential behavioral characteristics of the synovial fluid. We use the three models suggested by Hron et al. [ J. Hron, J. Málek, P. Pustějovská, K. R. Rajagopal, On concentration dependent shear-thinning behavior in modeling of synovial fluid flow, Adv. in Tribol.(In Press).] to study the problem, by appealing to a semi-inverse method. The assumed structure for the velocity field automatically satisfies the constraint of incompressibility, and the balance of linear momentum is solved together with a convection-diffusion equation. The results are compared to those associated with the Newtonian model. We also study the case in which an external pressure gradient is applied along the axis of the cylindrical annulus.

## 1. INTRODUCTION

There are several biological fluids as well as chemically-reacting polymeric liquids wherein the response characteristics of the fluids change due to the extent of the chemical reaction. In a biological material, such as blood, which in a homogenized sense is modeled as a fluid, the number of constituents are manifold: plasma, red and white blood cells, platelets, and a plethora of proteins; moreover, there are literally dozens of biochemical reactions that take place in order to maintain the blood in a delicate state of balance. For such complex fluids, the modeling advocated in this paper would be inappropriate. However, there are other biological liquids such as the synovial fluid which are less complex in structure in that there are fewer constituents and there are far fewer biochemical reactions. The present study addresses such biological fluids or polymeric liquids in which few or preferably one dominant chemical reaction takes place.

Synovial fluid is a mixture that responds in a viscoelastic fashion (see [1], [2], [3], and [4]). That the properties of the fluid change due to pathological conditions has been established by Gomez and Thurston [5]. When the synovial fluid is flowing under conditions where there are no instantaneous inputs, it essentially behaves as a generalized viscous fluid, namely a Stokesian fluid, and it is only when it is subject to instantaneous inputs does its viscoelastic nature manifest itself. Thus, the study of the flow of a generalized viscous fluid has some bearing on the flows of the synovial fluid. Recently, Hron et al. [6] studied the flow of three different fluid models which could be used as a model for the synovial fluid. All the models that they used fall into the category of generalized viscous fluids. One of them is a shear-thinning fluid model in which the power-law exponent is dependent on the concentration which leads to very interesting issues concerning its mathematical features. It is important to assess which of these three models studied by Hron et al. [6] is better suited to model the synovial fluid. We seek a better understanding of the models; especially with regard to determining if there are significant differences in the results between the models when the flow takes place in a totally different geometry. We consider an idealized geometry which allows us to use a semi-inverse method to study the problem while Hron et al. [6] use a finite element approach. However, unlike the geometry which Hron et al. [6] used, namely the flow in a rectangular region, here we consider a geometry wherein the boundary has a curvature which presents a new feature with regard to the flow domain in which to assess the characteristics of the model.

Synovial fluid is a mixture of plasma that is free of large proteins but containing polysaccharide molecules called Hyaluronan (HA). A detailed description of the composition and structure of the synovial fluid can be found in [7]. The presence of HA makes the fluid non-Newtonian, and the response of the synovial fluid changes with the concentration of HA. However, in view of the mass/volume fraction of the HA in the synovial fluid being small, we can approximate it essentially as a coexisting mixture of the plasma (which will henceforth be identified as the fluid) and a reactant (which is essentially HA). Let us consider a fluid whose property changes due to its reacting with a chemical agent whose mass fraction is significantly smaller than that of the fluid and furthermore the diffusion velocity of the chemical agent with respect to the fluid is far less than that of the fluid velocity. This implies that the relative velocity between the fluid and the diffusant is not significant, and the diffusant is essentially carried by the fluid. It is important to recognize that while the diffusion velocity with respect to the fluid is small, it cannot be ignored, for otherwise we would not have the diffusion taking place within the fluid. Essentially what we have, which is the counterpart to the diffusion of a fluid through a porous solid, is the diffusion of one fluid through another fluid that is flowing with a significantly higher velocity with respect to the relative velocity between the fluid and the diffusant. With respect to the velocity of the fluid, the diffusion velocity of the diffusant can be ignored. However, we do have a term that takes into account the changes in the concentration of the diffusant due to concentration gradients. The properties of the fluid change due to the concentration of the reactant, as the fluid chemically reacts with the diffusant. The study that is carried out is thus relevant not just to the flow of the synovial fluid, but also to polymeric fluids that are chemically-reacting with a diffusing substance, in our case (HA). Thus the approach used by Bridges et al. [8] can be used to study the flow in question, which is the methodology adopted by Hron et al. [6].

We study the flow of this reacting fluid mixture in the annular region between two coaxial cylinders, one of which is rotating; there is also flow due to the presence of a pressure gradient. The organization of the paper is as follows. The preliminaries and kinematics required for this paper are discussed first in section (2). The models chosen to characterize the response of the synovial fluid, and the corresponding governing equations are described in section (3) and section (4), respectively. In the following section, we describe the initial boundary value problem, the non-dimensionalization scheme, and the solution procedure in detail. The results are discussed in the final section.

## 2. KINEMATICS

Let $\mathcal{B}$ denote an abstract body, and let $\kappa_R$ and $\kappa_t$ represent the mappings, referred to as placers, which take the abstract body into the reference configuration $\kappa_R(\mathcal{B})$ and current configuration $\kappa_t(\mathcal{B})$, respectively. We shall choose to use the initial configuration as the reference configuration for the sake of convenience. Furthermore, we shall assume that there is a one-parameter family of placers which can be associated with the abstract motion of the body. We shall associate the parameter with time $t$, $t \in \mathbb{R}_0^+$. In view of this assumption, we can define a mapping $\boldsymbol{\chi}_{\kappa_R} : \kappa_R(\mathcal{B}) \times \mathbb{R}_0^+ \to \mathcal{E}$ such that:

$$\mathbf{x} = \boldsymbol{\chi}_{\kappa_R}(\mathbf{X}, t) \quad \text{where} \quad \mathbf{X} \in \kappa_R(\mathcal{B}) \quad \text{and} \quad \mathbf{x} \in \kappa_t(\mathcal{B}). \tag{1}$$

The velocity of the particle $\mathbf{X}$ is defined through the relation:

$$\mathbf{v}_{\kappa_R}(\mathbf{X}, t) = \frac{\partial \boldsymbol{\chi}_{\kappa_R}}{\partial t}(\mathbf{X}, t) := \frac{d\boldsymbol{\chi}_{\kappa_R}}{dt}, \tag{2}$$

where $\dfrac{d(\cdot)}{dt}$ represents the material time derivative, and the deformation gradient is defined as follows:

$$\mathbf{F}_{\kappa_R}(\mathbf{X}, t) = \frac{\partial}{\partial \mathbf{X}}[\boldsymbol{\chi}_{\kappa_R}(\mathbf{X}, t)] := \text{Grad}[\boldsymbol{\chi}_{\kappa_R}(\mathbf{X}, t)]. \tag{3}$$

We shall assume that $\det(\mathbf{F}_{\kappa_R}) > 0$. Furthermore, we assume that the motion is sufficiently smooth so that all of the derivative operations that follow have meaning. Since the motion is presumed to be one-to-one for each instant of time, the mapping (1) can be inverted so that:

$$\mathbf{X} = \boldsymbol{\chi}_{\kappa_R}^{-1}(\mathbf{x}, t), \tag{4}$$

and therefore any physical quantity of interest $\varsigma$ can be expressed in the various forms:

$$\varsigma = \tilde{\varsigma}(\mathbf{X}, t) = \bar{\varsigma}(\boldsymbol{\chi}_{\kappa_R}^{-1}(\mathbf{x}, t), t) = \hat{\varsigma}(\mathbf{x}, t). \tag{5}$$

We shall primarily be interested in stating the functions in terms of $\mathbf{x}$, i.e., we wish to study the mechanics of the reacting fluid mixture from an Eulerian perspective. Furthermore, we shall use lower case letters to represent differential operators such as $\text{div}(\cdot)$, $\text{grad}(\cdot)$, etc, which are carried out with respect to $\mathbf{x}$, and this is consistent with using capital letters to signify that the derivative is taken with respect to $\mathbf{X}$, as done previously with the definition of the deformation gradient (Lagrangean gradient of the motion). In addition, we shall henceforth omit the subscript $\kappa_R$.

It follows that

$$\frac{d}{dt}[\det(\mathbf{F})] = \det(\mathbf{F})\text{div}(\mathbf{v}), \tag{6}$$

and

$$\mathbf{L} := \text{grad}(\mathbf{v}) = \dot{\mathbf{F}}\mathbf{F}^{-1}, \tag{7}$$

where $\dot{(\ )}$ denotes the material time derivative. Now, if the body is incompressible and is therefore only capable of undergoing isochoric motions, then $\det(\mathbf{F}) = 1$ for each instant of time, and the relationship (6) reduces to

$$\text{div}(\mathbf{v}) = 0. \tag{8}$$

For the purpose of this work, we shall find it convenient to define the Rivlin-Ericksen tensors (see [9]). The first Rivlin-Ericksen tensor has the following representation:

$$\mathbf{A}_1 = \mathbf{L} + \mathbf{L}^T = 2\mathbf{D}, \tag{9}$$

where $\mathbf{D}$ is the symmetric part of the velocity gradient, i.e., the stretching tensor. The higher order Rivlin-Ericksen tensors are defined through the recurrence relationship:

$$\mathbf{A}_k = \frac{d\mathbf{A}_{k-1}}{dt} + \mathbf{A}_{k-1}\mathbf{L} + \mathbf{L}^T\mathbf{A}_{k-1}, \quad k \geq 2. \tag{10}$$

We shall employ a constrained mixture theory to study the behavior of synovial fluid which shear-thins or thickens due to the chemical reactions which take place. While the synovial fluid is comprised of numerous constituents: plasma, HA, fats, and leukocytes, the rheological properties of the mixture primarily depend upon the amount of plasma and HA. Thus, it is assumed that at each point within the mixture, the plasma (fluid) coexists with the HA (reactant), and the flow is constrained so these constituents move together. Henceforth, when referring to the "fluid" we shall mean the plasma, and we use the term "reactant" to stand for HA. We shall account for thinning and thickening by allowing the viscosity of the fluid to depend upon the shear-rate and the concentration, which changes by virtue of the chemical reactions that take place within the mixture. Here we define the concentration in the following manner:

$$c = \frac{\varrho_r}{\varrho_r + \varrho_f}, \tag{11}$$

where the quantities $\varrho_f$ and $\varrho_r$ represent the densities of the fluid and the coexisting reactant, respectively. We shall assume that the fluid is incompressible, and therefore only isochoric motions are possible.

## 3. Constitutive Assumptions

As mentioned in the introduction, the synovial fluid is a chemically-reacting, non-linearly viscoelastic, incompressible fluid mixture which possesses instantaneous elasticity. In fact, numerous experiments have been carried out which illustrate these facts ([1], [2], [4], [5], and [3]). However, if the fluid is not subject to instantaneous inputs, and if the inputs are sufficiently slow, then the fluid can be approximated as an isotropic and incompressible fluid of the differential type. The fluids of the differential type are incapable of instantaneous elastic response. The Cauchy stress for such a fluid has the following representation (see Truesdell and Noll [10]):

$$\mathbf{T} = -\pi\mathbf{1} + \mathbf{g}(\mathbf{A}_1, \mathbf{A}_2, ..., \mathbf{A}_k), \tag{12}$$

where $-\pi\mathbf{1}$ denotes the spherical stress due to the incompressibility constraint. Fluids that can be represented by this model are usually referred to as incompressible Rivlin-Ericksen fluids of complexity $k$, and these fluid models are frame-indifferent as the tensors $\mathbf{A}_k$ are frame-indifferent. For a more general class of fluids of the differential type of class $(m, k)$, see Dunn and Rajagopal [11]. From the standard representation theorems, for an isotropic, tensor-valued function it follows that the Cauchy stress for an incompressible Rivlin-Ericksen fluid of complexity 1 has the following form:

$$\mathbf{T} = -\pi\mathbf{1} + \phi_1\mathbf{A}_1 + \phi_2\mathbf{A}_1^2, \tag{13}$$

where $\phi_1$ and $\phi_2$ are each functions of the invariants $\text{tr}(\mathbf{A}_1^2)$ and $\text{tr}(\mathbf{A}_1^3)$. We are interested in a very special sub-class of (13) that is usually referred to as power-law fluids and such fluids can be used to describe shear-thinning a characteristic of the synovial fluid. For such fluids, the Cauchy stress reduces to:

$$\mathbf{T} = -\pi\mathbf{1} + \mu_a(\mathbf{A}_1)\mathbf{A}_1, \tag{14a}$$

and we shall assume that the apparent viscosity $\mu_a$ has the following structure:

$$\mu_a(\mathbf{A}_1) = \mu_0\left[1 + \gamma\mathrm{tr}(\mathbf{A}_1^2)\right]^n, \tag{14b}$$

where $\mu_0$, $\gamma$, $n$ are each constant. These fluids are capable of characterizing shear-thinning and shear-thickening phenomena, but cannot describe normal stress differences in simple shear flows. Furthermore, the model does not exhibit stress-relaxation effects nor does it possess the capability to describe instantaneous elasticity; the model is capable of characterizing certain types of non-linear creep, however. If the shear-index $n < 0$, then this model characterizes the behavior of a shear-thinning fluid while if $n > 0$, it characterizes the behavior of a shear-thickening fluid. If $n = 0$, then this model reduces to the incompressible, linearly viscous fluid model. The quantity $\mu_0$ is called the zero shear-rate viscosity and is defined as:

$$\mu_0 = \lim_{\kappa \to 0} \mu_a(\kappa), \tag{15}$$

where $\kappa$ is the shear-rate in a simple shear flow. The Lagrange multiplier $\pi$, in general, is a function of $\mathbf{x}$ and $t$, and since $\mathrm{tr}(\mathbf{A}_1) = 2\,\mathrm{div}(\mathbf{v}) = 0$, $\pi$ is equal to the mean normal stress. This is not the case with most incompressible, non-linear fluids as the trace of the extra stress is, in general, non-zero.

We shall assume that the fluid is initially of uniform temperature and is constrained to undergo only isothermal processes[*]. Consequently, the rate of dissipation $\Xi$ reduces to:

$$\Xi = \mathbf{T} \cdot \mathbf{D}, \tag{16}$$

where the scalar product between the Cauchy stress and the stretching tensor is often referred to as the stress-power. The $2^{nd}$ law of thermodynamics requires that the rate of dissipation, and in this case the stress-power, to be non-negative. If the stress-power is to be non-negative for all $n$, then we can conclude that $\gamma \geqq 0$. These fluids have been used in modeling the flows of polymeric liquids [12], and the mathematical properties pertaining to existence and uniqueness of the equations governing the flows of a fluid characterized by (14) have been analyzed in great detail (see [13]), and the stability of such flows has been explored by Málek et al. [14]. Finally, we note that such fluid models can also be developed within the thermodynamic framework recently proposed by Rajagopal and Srinivasa [15]. In fact, if it turns out that the temperature or viscoelasticity of the fluid should be taken into consideration, then the thermodynamic framework used in [15] is much more appropriate and leads to a much richer class of models.

As stated previously, we are interested in capturing the response of the synovial fluid mixture which not only shear-thins, but also chemically-thickens, and therefore we are interested in generalizations of the model (14). The model is generalized to account for chemical-thickening by allowing the apparent viscosity to also depend upon the concentration. There are a number of ways one could go about doing this; we shall be interested in three different models. The first model is due to Lai et al. [16] and Laurent et al. [17], and it allows the zero shear-rate viscosity

---

[*]The human body temperature continually fluctuates, and both the magnitude and manner in which these fluctuations take place can depend on a number of variables, both internal and external.

to change with concentration. In particular, the zero shear-rate viscosity takes the following form:

$$\mu_0(c) := \overline{\mu}_0 e^{\alpha c}, \tag{17}$$

where $\overline{\mu}_0$ is the viscosity of the plasma, $\alpha$ is a constant. The other models we shall consider allow for a much more intricate connection between the stretching and concentration, and these models are due to Hron et al. [6]. For these models, which we shall designate as model 2a and 2b, respectively, we allow the power-law exponent $n$ that appears in (14b) to depend on the concentration in the specific way summarized as follows:

$$n(c) := \frac{1}{2}\left(e^{-\alpha c} - 1\right), \tag{18}$$

and

$$n(c) := \sigma\left(\frac{1}{\alpha c^2 + 1} - 1\right), \tag{19}$$

respectively, where $\alpha$ and $\sigma$ are each constant. Also, for these models wherein the shear index depends upon the concentration, the zero shear-rate is constant, i.e., $\mu_0 = \overline{\mu}_0$. Hron et al. [6] show that the structure of these models allows for a better fit to relative viscosity $\left(\frac{\mu}{\mu_0}\right)$ / shear-rate data, for a variety of concentrations than the models discussed in [16] and [17]. Table 1 summarizes the constants used in each of the models. Thus, the general form of the models that we shall consider has the representation:

$$\mathbf{T} = -\pi\mathbf{1} + \mu(c, \mathbf{A}_1)\mathbf{A}_1, \tag{20a}$$

where

$$\mu(c, \mathbf{A}_1) := \mu_0(c)\left[\beta + \gamma \mathrm{tr}(\mathbf{A}_1^2)\right]^{n(c)}. \tag{20b}$$

Note that the constant $\beta$ has also been incorporated into the model (20) which allows for a better fit to the test data [6]. A sufficient condition for the rate of dissipation to be non-negative, for all $n(c)$, is to require both $\gamma \geqq 0$ and $\beta \geqq 0$, and the parameter values obtained from the regression analyses obviously agree with this constraint.

TABLE 1. Values for the Parameters used in each Model

| Parameter | Model 1 | Model 2a | Model 2b |
|---|---|---|---|
| $\alpha$ | 21.3 | 3.3 | 31.0 |
| $\beta$ | 1 | $7.1 \times 10^{-9}$ | $1.3 \times 10^{-8}$ |
| $\gamma$ [‡] | $\frac{6.96}{4}$ | $\frac{5.8 \times 10^{-8}}{4}$ | $\frac{8.5 \times 10^{-8}}{4}$ |
| $\sigma$ | – | – | 0.44 |
| $n$ | $-0.28$ | – | – |

[‡] Note that $\gamma$ defined in model (20) is $\frac{1}{4}$ times that used by Lai et al. [16], Laurent et al. [17], and Hron et al. [6] as their models are written in terms of $\mathbf{D}$ and the models herein are given in terms of $\mathbf{A}_1$.

## 4. Governing Equations

As stated previously, we employ the framework of constrained mixtures to model the synovial fluid mixture, and the constituents are therefore constrained to move together. Thus, we need not concern ourselves with the equations of motion for the reactant, just the equations which govern the flow of the fluid and the reactions which take place therein. We shall assume that the plasma is incompressible, and therefore we must ensure that the constraint (8) is satisfied. Thus, the balance of mass for the fluid:

$$\frac{d\varrho_f}{dt} + \varrho_f \text{div}(\mathbf{v}) = 0, \tag{21}$$

reduces to $\frac{d\varrho_f}{dt} = 0$. This implies that the fluid density is constant at each material point $\mathbf{X}$. However, since the fluid (plasma) is assumed to be homogeneous, it follows that $\varrho_f$ is the same constant everywhere. On substituting the constitutive relation (20) into the balance of linear momentum for the fluid:

$$\varrho_f \left[ \frac{\partial \mathbf{v}}{\partial t} + \text{grad}(\mathbf{v})\mathbf{v} \right] = \text{div}(\mathbf{T}) + \varrho_f \mathbf{b}, \tag{22}$$

we obtain:

$$\varrho_f \left[ \frac{\partial \mathbf{v}}{\partial t} + \text{div}(\mathbf{v} \otimes \mathbf{v}) \right] = -\text{grad}(\pi) + \text{div}[\mu(c, \mathbf{A}_1)\mathbf{A}_1]. \tag{23}$$

Note that we have utilized the condition (8) and neglected the body force $\varrho_f \mathbf{b}$ in deriving equation (23). We shall assume that the chemical reactions which take place within the mixture are governed by the convection-diffusion equation:

$$\frac{\partial c}{\partial t} + \text{div}(c\mathbf{v}) = -\text{div}(\mathbf{p}), \tag{24}$$

where $\mathbf{p}$ is a flux vector that is related to the reactions that take place within the synovial fluid. We shall assume that the flux vector $\mathbf{p}$ is given by a constitutive relation that is similar to that used in Fick's assumption, i.e.,

$$\mathbf{p} = -\mathbf{K}\,\text{grad}(c), \tag{25}$$

where $\mathbf{K}$ is, in general, a tensor-valued function of the concentration and the first Rivlin-Ericksen tensor. Therefore, the convection-diffusion equation reduces to:

$$\frac{\partial c}{\partial t} + \text{grad}(c) \cdot \mathbf{v} = \text{div}\left[\mathbf{K}(c, \mathbf{A}_1)\text{grad}(c)\right] \tag{26a}$$

after enforcing the constraint (8). In fact, the diffusivity could depend on $\mathbf{A}_k$, $k \geqq 2$, but we shall not consider such a possibility here. Experimental data suggests that, for the synovial fluid, the dependence on the concentration and the first Rivlin-Ericksen tensor is mild, and therefore we shall assume the diffusivity has the following form:

$$\mathbf{K}(c, \mathbf{A}_1) = D_c \mathbf{1}, \tag{26b}$$

where $D_c$ is a constant.

Thus, in general, the isochoricity constraint (8), the balance of linear momentum (23), and the convection-diffusion equation (26) are coupled, and one must solve them together for the velocity $\mathbf{v}$, the pressure $\pi$, and the concentration $c$ on a complicated, time-varying domain. Hron et al. [6] used a finite element approach to study the problem, and their approach can be used in complicated flow geometries. In this study, our interest lies in carrying out the analysis in a different but simpler geometry in view of the difficulties involved with solving the problem in a complex geometry, especially one in which the curvature effects come into play. Thus, we shall consider a relatively simple initial boundary value problem.

## 5. An Initial Boundary Value Problem

5.1. **Assumptions for the Velocity Field.** As with any semi-inverse approach, one makes an attempt to reduce the complexity of the governing equations so that they become more tractable in a coordinate system appropriate for the problem of interest. The aim is to glean some insight into the full set of equations by studying the solution of the resulting simplified equations, usually on an infinite domain. For the problem under consideration, one such approach would be to study the plane flow wherein the velocity field has the form:

$$\mathbf{v} = v(\nu, \xi, t)\mathbf{e}_\nu + w(\nu, \xi, t)\mathbf{e}_\xi, \tag{27}$$

in a confocal-ellipsoidal coordinate system $(\zeta, \nu, \xi)$. Of course, in such a coordinate system the lines of constant $\zeta$ form ellipsoids, while the lines of constant $\nu$ and $\xi$ form hyperboloid of one sheet and two sheets, respectively; these surfaces could characterize the articular membrane which surrounds the bone and the synovial membrane of the joint. The drawback of such an approach is that as the appendage moves, the bone rotates and the coordinate system is no longer useful as the boundaries do not coincide with lines of constant hyperboloids.

Another approach is to study the response of the reacting fluid mixture that fills to the gap between two spheres and allow the outer sphere to rotate while the inner sphere remains fixed. One could utilize the similarity transformation suggested by Bandelli and Rajagopal [18]:

$$\mathbf{v} = v(r, t)\sin(\theta)\mathbf{e}_\varphi, \tag{28}$$

where $(r, \theta, \varphi)$ is a spherical coordinate system. Such an assumption can only be made if the unsteady part in the inertial term is neglected, or the inertial term is neglected as a whole.

We consider the unsteady flow of the reacting mixture which is confined to the space between two cylinders, i.e., the mixture domain $\Omega$ is defined as: $r_i \leqq r \leqq r_o$, $0 \leqq \theta < 2\pi$, and $-\infty < z < \infty$ (see Figure 1). Here $(r, \theta, z)$ represents a cylindrical-polar coordinate system that is placed such that $r = r_i$ is the boundary of the inner cylinder, which is fixed, and $r = r_o$ is the inner wall of the outer cylinder (annulus), which is free to rotate. We shall employ a semi-inverse approach and seek a velocity field of the form:

$$\mathbf{v} = v_\theta(r, t)\mathbf{e}_\theta + v_z(r, t)\mathbf{e}_z, \tag{29}$$

where $\mathbf{e}_\theta$ and $\mathbf{e}_z$ denote the base vectors in the $\theta$ and $z$ coordinate directions, respectively. Note that this assumption for the velocity field automatically satisfies the isochoricity constraint (8). On substituting the similarity transformation (29) into the balance of linear momentum (23), expressed in a cylindrical-polar coordinate system, we obtain:

$$\frac{\partial \pi}{\partial r} = \frac{\varrho_f}{r} v_\theta^2, \tag{30a}$$

$$\frac{1}{r^2} \frac{\partial}{\partial r} \left[ r^2 \mu \left( c, \frac{\partial v_\theta}{\partial r}, \frac{\partial v_z}{\partial r} \right) \left( \frac{\partial v_\theta}{\partial r} - \frac{v_\theta}{r} \right) \right] = \varrho_f \frac{\partial v_\theta}{\partial t}, \tag{30b}$$

and

$$\frac{1}{r}\frac{\partial}{\partial r}\left[r\mu\left(c,\frac{\partial v_\theta}{\partial r},\frac{\partial v_z}{\partial r}\right)\frac{\partial v_z}{\partial r}\right] - \frac{\partial \pi}{\partial z} = \varrho_f \frac{\partial v_z}{\partial t}, \tag{30c}$$

where

$$\mu\left(c,\frac{\partial v_\theta}{\partial r},\frac{\partial v_z}{\partial r}\right) := \mu_0(c)\left\{\beta + 2\gamma\left[\left(\frac{\partial v_\theta}{\partial r} - \frac{v_\theta}{r}\right)^2 + \left(\frac{\partial v_z}{\partial r}\right)^2\right]\right\}^{n(c)}. \tag{30d}$$

Note that in deriving (30), we have tacitly assumed that the pressure is independent of the $\theta$ coordinate. It would then be natural to assume, since we shall consider a fully-developed, oscillating flow, an axial pressure gradient of the form:

$$\frac{\partial \pi}{\partial z} = -[A + B\cos(f_z t)], \tag{31}$$

where $A$ and $B$ are constants, and $f_z$ is the frequency of the oscillation. On integrating equation (31) we find the pressure:

$$\pi(r,z,t) = -[A + B\cos(f_z t)]z + h(r,t), \tag{32}$$

where $h(r,t)$ satisfies:

$$\frac{\partial h}{\partial r} = \frac{\varrho_f}{r}v_\theta^2, \tag{33}$$

in view of the radial component of the balance of linear momentum (30a). If the circumferential component of velocity is known, then $h(r,t)$ can be calculated via the formula:

$$h(r,t) = \varrho_f\int_{r_i}^r \frac{v_\theta^2(s,t)}{s}ds + a(t), \tag{34}$$

where $a(t)$ is an arbitrary function of time. The remainder of this paper is focused on determining $v_\theta$, $v_z$, and $c$. To do so, we require both initial conditions and boundary conditions for $v_\theta$, $v_z$, and $c$; let us first consider the conditions on the velocity.

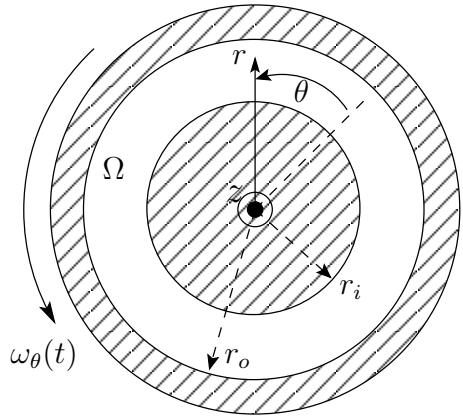

FIGURE 1. Illustration of the geometry for the problem under consideration along with a cylindrical-polar coordinate system.

5.2. **Initial Conditions and Boundary Conditions on the Velocity Field.** We shall assume that the fluid is initially at rest, and therefore the non-zero components of the velocity field must satisfy:

$$v_\theta(r, 0) = v_z(r, 0) = 0 \quad \forall \ r \in [r_i, r_o]. \tag{35}$$

For fluids such as water flowing under "normal conditions," it is customary to enforce the "no-slip" condition[§] at the boundaries. However, if the mixture is sufficiently dense in either plasma or HA, then generalizations of Navier's slip condition may be more appropriate. For example, following the approach recommended by Rajagopal [19], one could have a boundary condition of the form:

$$\mathbf{v}_{rel} \cdot \boldsymbol{\tau} + M(|\mathbf{T}^T \mathbf{n} \cdot \boldsymbol{\tau}|, c) \frac{\mathbf{T}^T \mathbf{n} \cdot \boldsymbol{\tau}}{|\mathbf{T}^T \mathbf{n} \cdot \boldsymbol{\tau}|} = 0, \tag{36}$$

where $\mathbf{v}_{rel}$ is the relative velocity between the boundary surface and the fluid, and $\boldsymbol{\tau}$ and $\mathbf{n}$ are the unit tangent and unit normal vectors, respectively, to the boundary surface. Note that $M$ has been modified to account for the concentration of HA. In addition to the condition (36), the usual condition $\mathbf{v} \cdot \mathbf{n} = 0$ should also be satisfied if one assumes that the velocity of the mass diffusing through the articular cartilage and synovial membrane is negligible with respect to the velocity of the mass in the cavity. One example of the function $M$ that has relevance for the problem under consideration is:

$$M(|\mathbf{T}^T \mathbf{n} \cdot \boldsymbol{\tau}|, c) = \frac{1 - \lambda(c)}{\lambda(c)} |\mathbf{T}^T \mathbf{n} \cdot \boldsymbol{\tau}|, \tag{37}$$

and one could require $\lambda(0)$ to be small so that when the fluid mixture primarily consists of plasma, the fluid nearly slips at the boundary. Also, if one requires $\lambda(1) = 1$, then when the mixture is dense with HA, the no-slip boundary condition is enforced. With that being said, if the mixture is dense enough with HA, then the model represented by equation (20) may be ineffective in describing the response of the synovial fluid as viscoelastic effects could play a more dominant role. As a first order approximation we shall choose to enforce the no-slip boundary condition at the walls by requiring:

$$v_\theta(r_i, t) = 0 \quad \text{and} \quad v_\theta(r_o, t) = r_o \omega_\theta(t), \quad t \geqq 0, \tag{38}$$

where $\omega_\theta$ is the angular velocity of the outer cylinder. We shall assume that $\omega_\theta$ takes the form: $\omega_\theta(t) = \overline{\omega}_\theta[1 - \cos(f_\theta t)]$, which could represent some cyclic routine such as exercise. Also, the cylinders are not permitted to translate, and therefore we must enforce the following boundary conditions for the axial component of the velocity:

$$v_z(r_i, t) = v_z(r_o, t) = 0, \quad t \geqq 0. \tag{39}$$

The solution to equations (30b) and (30c) ultimately depends upon the concentration, and this is the issue we take up next.

---

[§]While this condition is supposed to have had the backing of Stokes, he was far from unequivocal about its applicability to all flows of fluid such as water; he only advocated the same for sufficiently slow flows (see Rajagopal [19] for a detailed discussion of the same).

5.3. **Assumptions for the Concentration.** Since the mixture is confined to the annulur region between the two cylinders, and one of the cylinders is oscillating, we assume the concentration has the following structure:

$$c = c(r, t), \tag{40}$$

and therefore the governing equation (26) reduces to:

$$\frac{D_c}{r}\frac{\partial}{\partial r}\left(r\frac{\partial c}{\partial r}\right) = \frac{\partial c}{\partial t}. \tag{41}$$

5.4. **Initial Conditions and Boundary Conditions for the Concentration.** We shall require the initial concentration be constant $c_i$ over the entire domain. The synovial membrane lines the joint cavity with the exception of the joint itself which is surrounded by articular cartilage. If we take the inner cylinder to represent the joint, and the outer cylinder to be the cavity wall lined with synovial membrane, then the concentration must satisfy the following equations:

$$\frac{\partial c}{\partial r}(r_i, t) = 0 \;\; \text{and} \;\; c(r_o, t) = c_{r_o}(t), \;\; t \geqq 0. \tag{42}$$

Condition (42) presumes that there is no diffusion of HA through the articular cartilage, which does in fact act as an HA reservoir. A more appropriate condition at the outer boundary would be to prescribe the concentration, until the mean value of HA within a region near the outer boundary reaches an optimum value at which time no more HA is secreted. This boundary condition would take the form:

$$\left.\begin{array}{c} c(r_o, t) \\ \frac{\partial c}{\partial r}(r_o, t) \end{array}\right\} = \left\{\begin{array}{ll} \tilde{c}, & c_{mean} < \overline{c} \\ 0, & c_{mean} \geqq \overline{c} \end{array}\right., \tag{43a}$$

where

$$c_{mean}(t) := \frac{1}{r_o - \overline{r}}\int_{\overline{r}}^{r_o} c(r, t)dr, \tag{43b}$$

$\tilde{c}$ is a constant, and $\overline{c}$ represents the optimal concentration level of HA within the layer $r_o - \overline{r}$ adjacent to the synovial membrane. Obviously the boundary condition (43) is a complicated condition to satisfy numerically as one must know, at a given instant of time, the entire concentration profile. Thus, in addition to the iterations required due to the non-linearity of the equations, one must develop a iterative stencil to ensure that the condition (43) is satisfied, even if approximately. Thus, we shall first consider the more simple case wherein the conditions given by (42) must hold.

5.5. **Non-dimensionalization.** We now proceed to non-dimensionalize the governing equations by introducing the mapping $(r, z, t) \rightarrowtail (\hat{r}, \hat{z}, \hat{t})$ and parameters:

$$\hat{r} = \frac{r - r_i}{r_o - r_i}, \qquad \hat{z} = \frac{z}{L}, \qquad \hat{t} = f_\theta t, \qquad \hat{v} = \frac{v_\theta}{r_o \overline{\omega}_\theta}, \qquad \hat{w} = \frac{v_z}{r_o \overline{\omega}_\theta}, \qquad \hat{\pi} = \frac{\pi}{\varrho_f r_o \overline{\omega}_\theta L f_\theta},$$

$$c = c(r, t) = \hat{c}(\hat{r}, \hat{t}), \qquad \hat{\mu}_0 = \frac{\mu_0}{\overline{\mu}_0}, \qquad \hat{\omega}_\theta = \frac{\omega_\theta}{\overline{\omega}_\theta}, \qquad p_f = \frac{f_z}{f_\theta}, \qquad p_g = \frac{r_i}{r_o - r_i}, \tag{44}$$

$$p_\gamma = \frac{2\gamma r_o^2 \overline{\omega}_\theta^2}{(r_o - r_i)^2}, \qquad p_\beta = \beta,$$

where $L$ is a characteristic length. First, we record the non-dimensional solution for the pressure for completeness:

$$\hat{\pi}(\hat{r}, \hat{z}, \hat{t}) = -[p_A + p_B \cos(p_f \hat{t})]\hat{z} + \int_0^{\hat{r}} \frac{\hat{v}^2(\hat{s}, \hat{t})}{\hat{s} + p_g} d\hat{s} + \hat{a}(\hat{t}), \tag{45}$$

where $\hat{a}(\hat{t}) := \frac{a(t)}{\varrho_f r_o \overline{\omega}_\theta L f_\theta}$, and the pressure gradient coefficients are given by the formulae:

$$p_A := \frac{A}{\varrho_f r_o \overline{\omega}_\theta f_\theta} \quad \text{and} \quad p_B := \frac{B}{\varrho_f r_o \overline{\omega}_\theta f_\theta}. \tag{46}$$

The non-dimensional forms of the governing equations (30b), (30c), and (41) are:

$$\frac{1}{Re} \frac{1}{(\hat{r} + p_g)^2} \frac{\partial}{\partial \hat{r}} \left[ (\hat{r} + p_g)^2 \hat{\mu} \left( \hat{c}, \frac{\partial \hat{v}}{\partial \hat{r}}, \frac{\partial \hat{w}}{\partial \hat{r}} \right) \left( \frac{\partial \hat{v}}{\partial \hat{r}} - \frac{\hat{v}}{\hat{r} + p_g} \right) \right] = \frac{\partial \hat{v}}{\partial \hat{t}}, \tag{47a}$$

and

$$\frac{1}{Re} \frac{1}{(\hat{r} + p_g)} \frac{\partial}{\partial \hat{r}} \left[ (\hat{r} + p_g) \hat{\mu} \left( \hat{c}, \frac{\partial \hat{v}}{\partial \hat{r}}, \frac{\partial \hat{w}}{\partial \hat{r}} \right) \frac{\partial \hat{w}}{\partial \hat{r}} \right] + p_A + p_B \cos(p_f \hat{t}) = \frac{\partial \hat{w}}{\partial \hat{t}}, \tag{47b}$$

$$\frac{1}{Pe} \frac{1}{\hat{r} + p_g} \frac{\partial}{\partial \hat{r}} \left[ (\hat{r} + p_g) \frac{\partial \hat{c}}{\partial \hat{r}} \right] = \frac{\partial \hat{c}}{\partial \hat{t}}, \tag{47c}$$

respectively, where

$$\hat{\mu} \left( \hat{c}, \frac{\partial \hat{v}}{\partial \hat{r}}, \frac{\partial \hat{w}}{\partial \hat{r}} \right) := \hat{\mu}_0(\hat{c}) \left\{ p_\beta + p_\gamma \left[ \left( \frac{\partial \hat{v}}{\partial \hat{r}} - \frac{\hat{v}}{\hat{r} + p_g} \right)^2 + \left( \frac{\partial \hat{w}}{\partial \hat{r}} \right)^2 \right] \right\}^{n(\hat{c})}, \tag{47d}$$

and the Reynold's and Peclet numbers are defined as follows:

$$Re := \frac{\varrho_f f_\theta (r_o - r_i)^2}{\overline{\mu}_0} \quad \text{and} \quad Pe := \frac{f_\theta (r_o - r_i)^2}{D_c}, \tag{47e}$$

respectively. Below we summarize both the initial conditions and boundary conditions, in non-dimensional form, for $\hat{v}$ and $\hat{w}$:

$$\hat{v}(\hat{r}, 0) = \hat{w}(\hat{r}, 0) = 0, \quad \forall \ \hat{r} \in [0, 1], \tag{48a}$$

$$\hat{v}(0, \hat{t}) = 0 \quad \text{and} \quad \hat{v}(1, \hat{t}) = 1 - \cos(\hat{t}), \quad \hat{t} \geqq 0, \tag{48b}$$

$$\hat{w}(0, \hat{t}) = \hat{w}(1, \hat{t}) = 0, \quad \hat{t} \geqq 0, \tag{48c}$$

and for the concentration $\hat{c}$:

$$\hat{c}(\hat{r}, 0) = 0.1, \quad \forall \ \hat{r} \in [0, 1], \tag{49a}$$

$$\frac{\partial \hat{c}}{\partial \hat{r}}(0, \hat{t}) = 0, \quad \hat{t} \geqq 0 \quad \text{and} \quad \hat{c}(1, \hat{t}) = \begin{cases} 0.1 + 0.1\hat{t}, & \hat{t} \in [0, 2] \\ 0.3, & \hat{t} > 2 \end{cases}. \tag{49b}$$

Furthermore, we shall require $p_A = -p_B$, and $p_B = 1$, and these parameter values agree with the initial condition (48a). The ramp function given in condition (49b) is the same used by Hron et al. [6] and this, to some extent, represents physiological data.

## 6. Discussion

Equations (47a)–(47c) constitute a non-linear system of partial differential equations that must be solved subject to the previously discussed initial conditions and boundary conditions, for both the velocity and the concentration. We use the subroutine `pdepe` within the software package MATLAB. This package utilizes a Petrov-Galerkin finite element method to generate a system of ordinary differential equations [20]. The ordinary differential equations are then solved using standard methods.

The set of coupled partial differential equations are solved for three shear-thinning, chemically-thickening models (Model 1, Model 2a, and Model 2b), and the results are compared with those for the Newtonian model. Note that since the diffusivity was assumed to be constant, the concentration profiles are the same for all the models (see Figure (3d)), i.e., while the components of velocity depend upon the concentration, the converse is not true. This is because the equation governing concentration is not influenced by the velocity field as the diffusivity is a constant, and a special form is assumed for the concentration and the velocity field. If one allows for the diffusivity to depend on the norm of the first Rivlin-Ericksen tensor, $|\mathbf{A}_1|$, then the equations would be further coupled, and the concentration profiles would vary depending on the model, at any given instant of time. The concentration profiles will also further change if an external pressure gradient is applied in the $z$ direction, when one assumes the concentration to depend on the norm of the first Rivlin-Ericksen tensor. Furthermore, note that one cycle is complete in (non-dimensional) time of $2\pi$, regardless of whether or not there is an axial pressure gradient present.

We shall first consider the case wherein there is no pressure gradient applied in the $z$ direction. The non-dimensional velocity profiles in the $\theta$ direction are depicted in Figures (2), (3a), and (3b). The non-dimensional velocity in the $z$ direction is zero for all cycles for each of the models studied because there is no pressure gradient in the $z$ direction, and the axial velocity is zero at the boundaries. The non-dimensional apparent viscosity variations are compared in Figures (4) and (5). Obviously the plots for the chemically-reacting fluid illustrate a great departure from those associated with the Newtonian model, and this is to be expected in view of the change in the concentration over the domain (Figure (3d)). Now, as the number of cycles increases, more reactant enters the domain from the outer boundary, and this reactant stays in the domain as the inner cylinder is non-porous. This, of course, increases the concentration of the reactant in the region between the cylinders (see Figure (3d)), which in turn increases the apparent viscosity as observed in Figures (4) and (5). Hence, as the cycles continue, the region near the outer boundary becomes more and more "viscous." This outer layer of highly viscous fluid squeezes the less viscous fluid particles, radially inward, and causes them to accelerate in the $\theta$ direction. Such behavior is illustrated in Figures (2), (3a), and (3b), wherein the velocity in the $\theta$ direction increases in the central region with number of cycles. Such phenomena was also observed by Bridges and Rajagopal [8] while studying shear-thinning/chemical-thickening fluids. However, for their model, the zero shear-rate viscosity had a quadratic dependence on the concentration, and the shear-index was assumed constant. Furthermore, it was assumed that the diffusivity depends on the first Rivlin-Ericksen tensor. In addition, the geometry and boundary conditions were slightly different than those for the problem under consideration.

Next, we shall study the case when there is an external pressure gradient in the $z$ direction. The non-dimensional velocity profiles in the $\theta$ direction are shown in Figures (6) and (7), and the non-dimensional velocity profiles in the $z$ direction are compared in Figures (8), (9a), and (9b) for the various models under consideration. The non-dimensional apparent viscosities in this case are compared in Figures (10) and (11). One can also infer from Figures (8), (9a), (9b), and (9d), that as the number of cycles increase, the component of velocity in the $z$ direction decreases in the central region. Now, $\hat{w}$ is zero at the boundaries, and is driven solely by the pressure gradient. As the number of cycles increases, there is a significant increase in the apparent viscosity in the central region. This increase in viscosity greatly mitigates the driving force of the pressure gradient, and thus decreases the velocity $\hat{w}$. Furthermore, in comparison with the results from the case when there is no pressure gradient, it can be seen that at fewer cycles, there is significant difference in the velocity profiles in the $\theta$ direction and the apparent viscosity. For instance, one can see this difference by comparing Figures (2d) and (6d) at 3.5 cycles. However, for more cycles (12.5, 34.5), the plots look similar. The reason being that the velocity in the $z$ direction, and its gradients are steep, at lower cycles. This phenomena affects the apparent viscosity, and hence in turn the velocity in the $\theta$ direction. At higher cycles, as discussed above the velocity in the $z$ direction is small, and its gradients are small as well; thus, there is no significant effect on the apparent viscosity and the velocity in the $\theta$ direction.

When one compares the three shear-thinning, chemically-thickening models, the velocity and the apparent viscosity profiles for Model 2a and Model 2b are similar, and those of Model 1 are quite distinct from Models 2a, 2b; this was something Hron et al. [6] also observed. It is also interesting to note that the nature of the shapes for the velocity plots for shear-thinning, chemically-thickening models, and the Newtonian model is completely different (see Figures (2), (3a), and (3b), for instance) - the plots for Models 1, 2a, 2b are concave and the plots for the Newtonian model are convex. It is quite possible that this is due to the curvature in the flow domain since the results for a rectangular domain as seen in [6] do not show such a behavior. In the study by Hron et al. [6], they find that the structure of the solutions for all the solutions for the three models are similar. Thus, studying the fluid in a different geometry seems to provide some insight in that the structure of the solution for the non-Newtonian fluids are markedly different from those of the Newtonian fluid. We were hoping that the curvature effects would help in providing us with new information concerning the response of these non-Newtonian fluid models and our guess was correct. These non-Newtonian fluids are quite different from the Newtonian fluid as evidenced by the profiles in the non-Newtonian case being concave and the Newtonian case being convex. This implies that the nature of the shear stresses in the non-Newtonian case and the Newtonian case are totally different. It would be useful to find an initial boundary value problem where the three different non-Newtonian models lead to different solutions with the hope that one agrees better with the corresponding experimental results than the two others, thereby allowing us to choose an appropriate model.

This study is merely a preliminary study to assess a fluid model that can describe the response of fluids such as the synovial fluid. This study will be followed by an analysis of such a fluid in a realistic geometry, the resolution of which will involve numerical resolution based on appealing to methods such as the finite element method.

CRAIG BRIDGES, TEXAS A&M UNIVERSITY, DEPARTMENT OF MECHANICAL ENGINEERING, 3123 TAMU, COLLEGE STATION TX 77843-3123, UNITED STATES OF AMERICA
    *E-mail address*: craig_bridges@tamu.edu

SATISH KARRA, TEXAS A&M UNIVERSITY, DEPARTMENT OF MECHANICAL ENGINEERING, 3123 TAMU, COLLEGE STATION TX 77843-3123, UNITED STATES OF AMERICA
    *E-mail address*: satkarra@tamu.edu

K. R. RAJAGOPAL (CORRESPONDING AUTHOR), TEXAS A&M UNIVERSITY, DEPARTMENT OF MECHANICAL ENGINEERING, 3123 TAMU, COLLEGE STATION TX 77843-3123, UNITED STATES OF AMERICA, PHONE: 1-979-862-4552, FAX: 1-979-845-3081
    *E-mail address*: krajagopal@tamu.edu

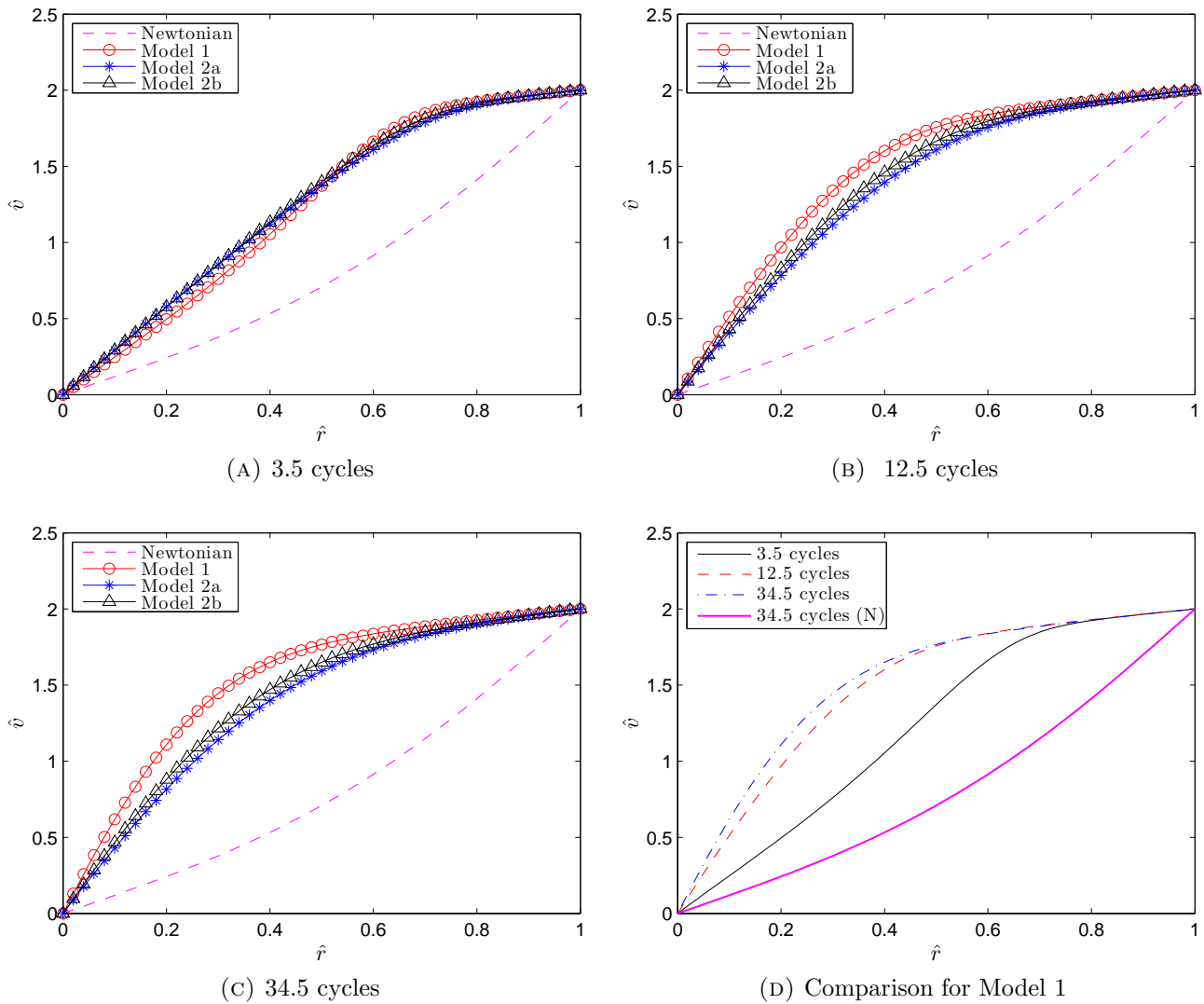

(A) 3.5 cycles

(B) 12.5 cycles

(C) 34.5 cycles

(D) Comparison for Model 1

FIGURE 2. Non-dimensional velocity profiles in the $\theta$ direction. Comparison for different models at different cycles are shown in (a), (b), (c). In (d), the non-dimensional velocity profiles in the $\theta$ direction are compared for different cycles for Model 1. The parameters chosen are $r_i = 1$, $r_o = 1.2$, $p_A = -p_B = 0$, $p_f = 1$, $\overline{\omega}_\theta = 1$, $Re = 10$, $Pe = 1000$. The pressure gradient in the $z$ direction is chosen to be zero here.

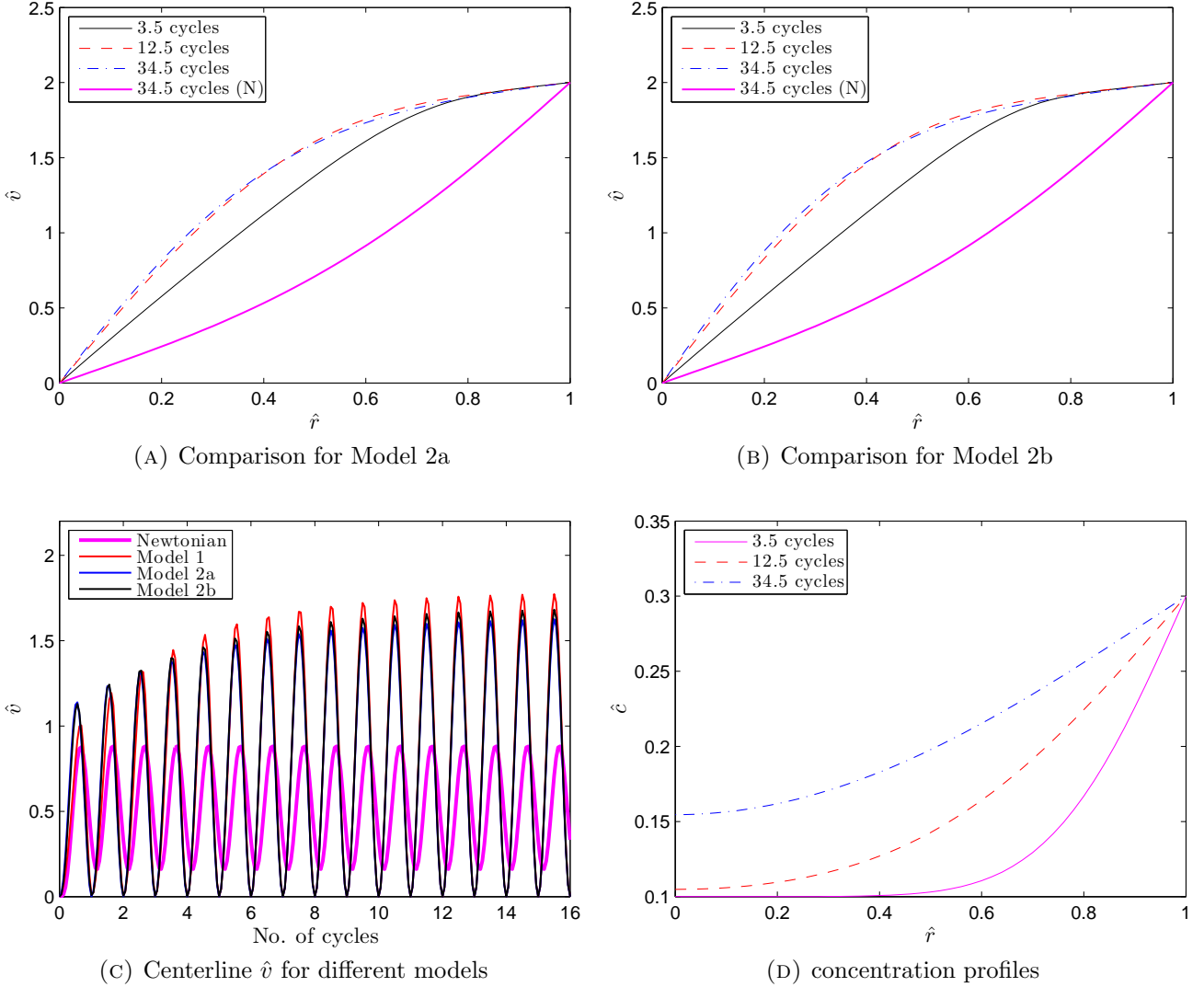

(A) Comparison for Model 2a

(B) Comparison for Model 2b

(C) Centerline $\hat{v}$ for different models

(D) concentration profiles

FIGURE 3. The non-dimensional velocity profiles in the $\theta$ direction are compared for different cycles. Results for Model 2a and Model 2b are shown in (a) and (b) respectively. Non-dimensional centerline velocity (at $\hat{r} = 0.5$) in the $\theta$ direction as a function of number of cycles for different models is shown in (c). Concentration profiles at different cycles are shown in (d). Since, the diffusivity is chosen to be a constant, the concentration profiles are the same for all the models. The parameters chosen are $r_i = 1$, $r_o = 1.2$, $p_A = -p_B = 0$, $p_f = 1$, $\overline{\omega}_\theta = 1$, $Re = 10$, $Pe = 1000$. The pressure gradient in the $z$ direction is chosen to be zero here.

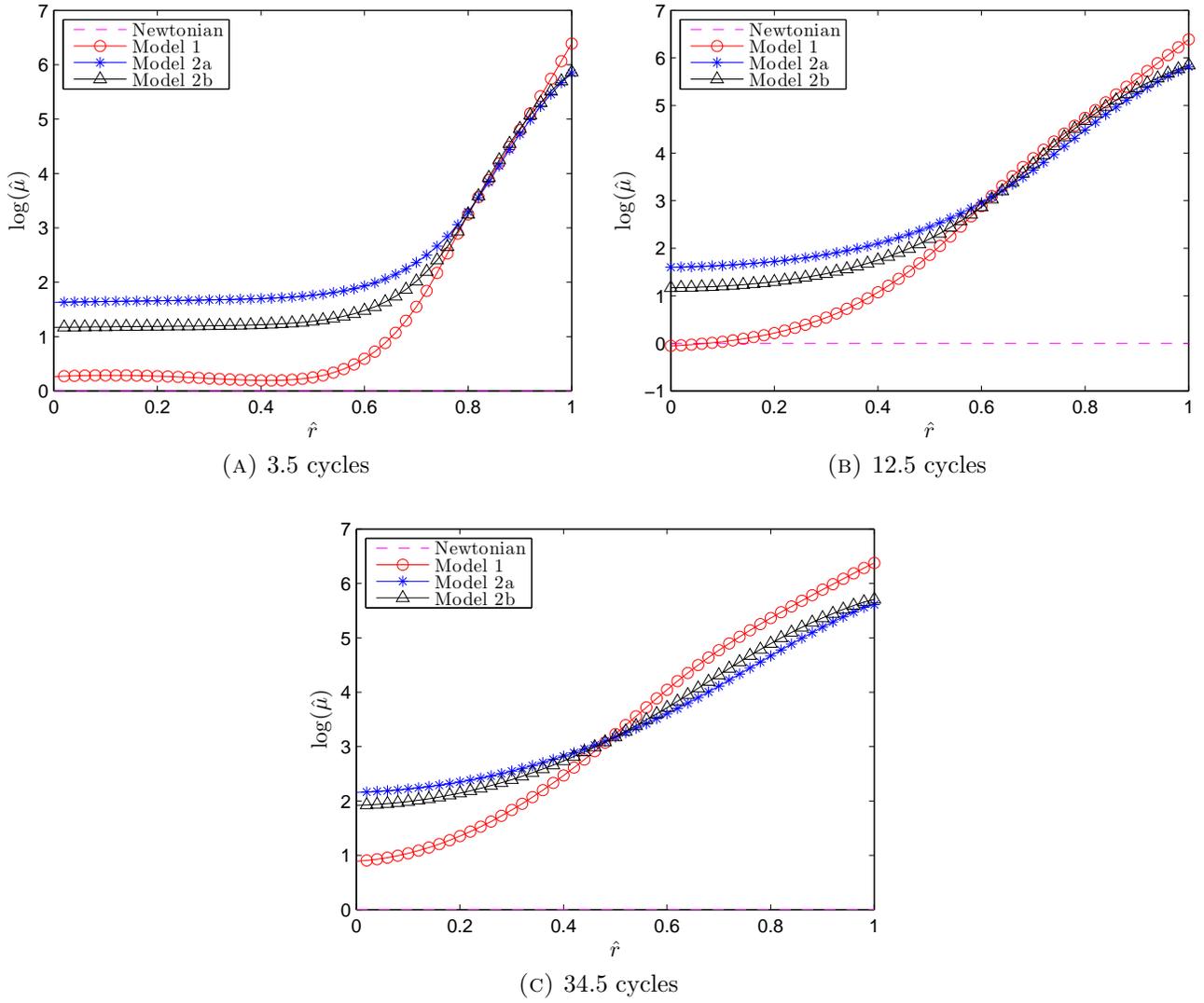

(A) 3.5 cycles

(B) 12.5 cycles

(C) 34.5 cycles

FIGURE 4. Non-dimensional apparent viscosity profile comparison for different models at different cycles are shown in (a), (b), (c). The parameters chosen are $r_i = 1$, $r_o = 1.2$, $p_A = -p_B = 0$, $p_f = 1$, $\overline{\omega}_\theta = 1$, $Re = 10$, $Pe = 1000$. The pressure gradient in the $z$ direction is chosen to be zero here.

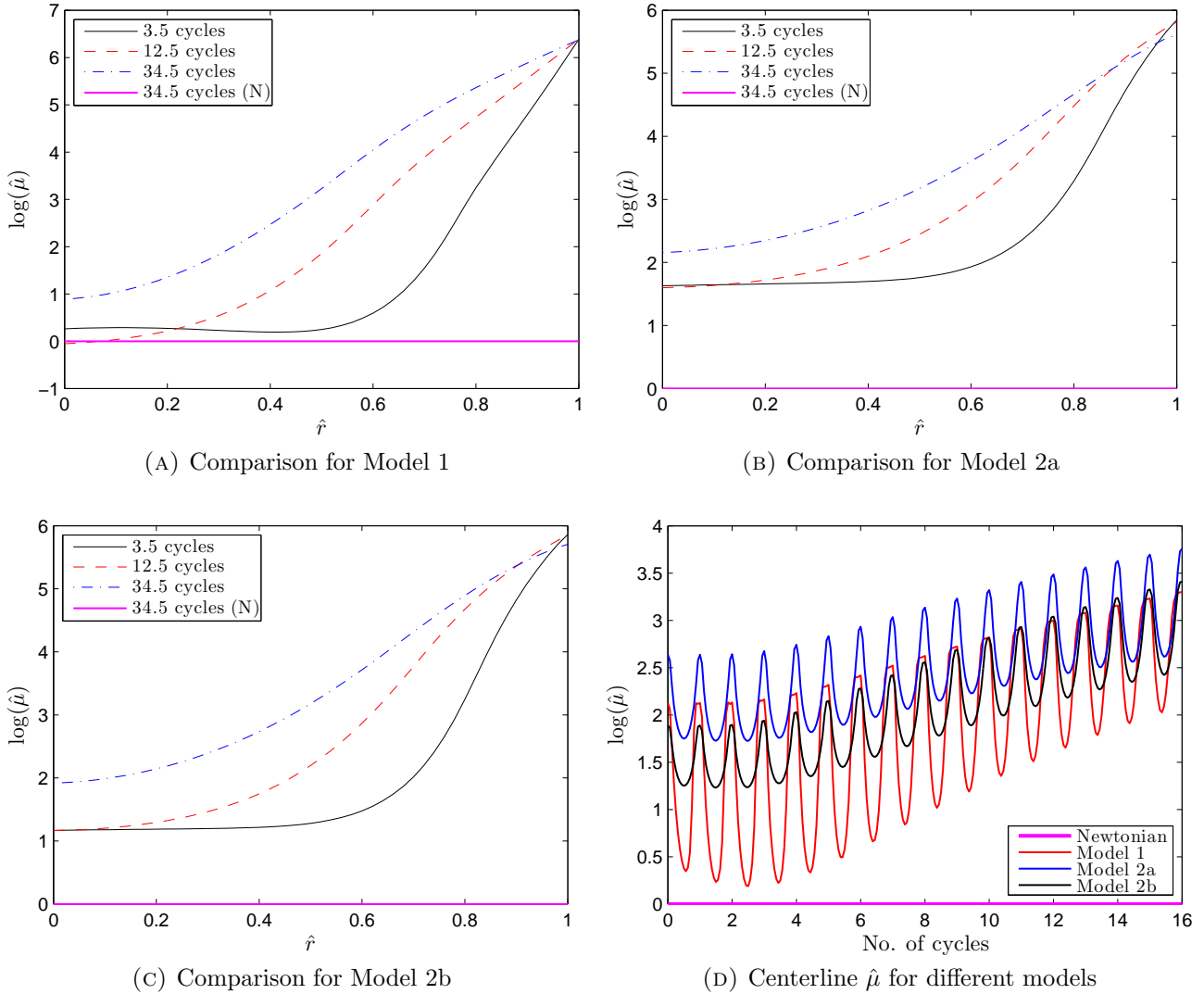

(A) Comparison for Model 1

(B) Comparison for Model 2a

(C) Comparison for Model 2b

(D) Centerline $\hat{\mu}$ for different models

FIGURE 5. The non-dimensional apparent viscosity profiles are compared for different cycles for Model 1, Model 2a and Model 2b. Results are shown in (a), (b) and (c) respectively. Non-dimensional centerline apparent viscosity (at $\hat{r} = 0.5$) as a function of number of cycles for different models is shown in (d). The parameters chosen are $r_i = 1$, $r_o = 1.2$, $p_A = -p_B = 0$, $p_f = 1$, $\overline{\omega}_\theta = 1$, $Re = 10$, $Pe = 1000$. The pressure gradient in the $z$ direction is chosen to be zero here.

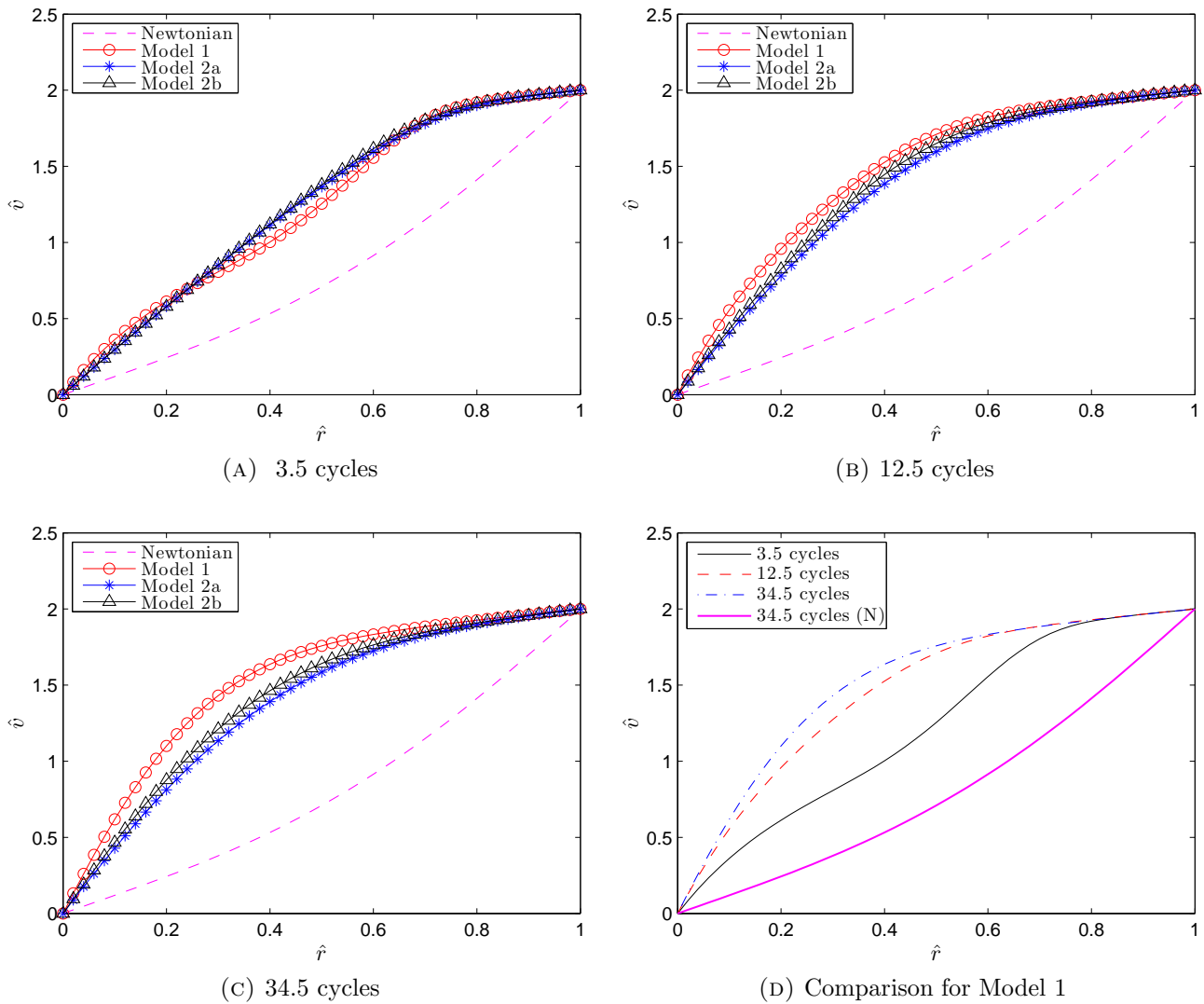

(A)  3.5 cycles

(B)  12.5 cycles

(C)  34.5 cycles

(D)  Comparison for Model 1

FIGURE 6.  Non-dimensional velocity profiles in the $\theta$ direction. Comparison for different models at different cycles are shown in (a), (b), (c). In (d), the non-dimensional velocity profiles in the $\theta$ direction are compared for different cycles for Model 1. The parameters chosen are $r_i = 1$, $r_o = 1.2$, $p_A = -p_B = 1$, $p_f = 1$, $\overline{\omega}_\theta = 1$, $Re = 10$, $Pe = 1000$. The concentation profiles for this set of results is same as (3d) because the diffusivity is chosen to be constant.

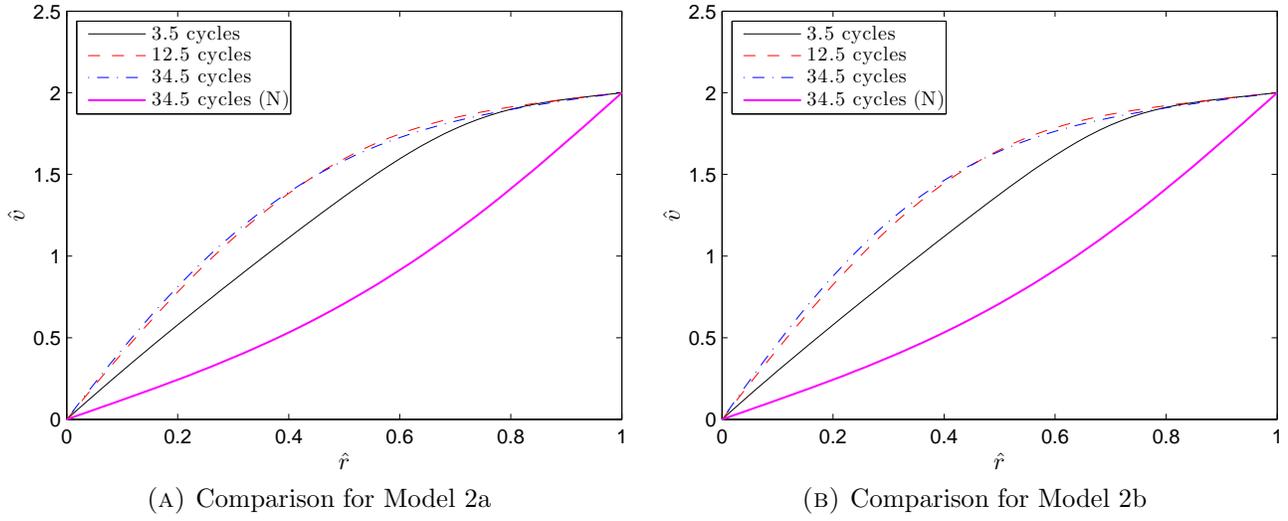

(A) Comparison for Model 2a                    (B) Comparison for Model 2b

FIGURE 7. The non-dimensional velocity profiles in the $\theta$ direction are compared for different cycles for Model 2a and Model 2b are shown in (a) and (b) respectively. The parameters chosen are $r_i = 1$, $r_o = 1.2$, $p_A = -p_B = 1$, $p_f = 1$, $\overline{\omega}_\theta = 1$, $Re = 10$, $Pe = 1000$.

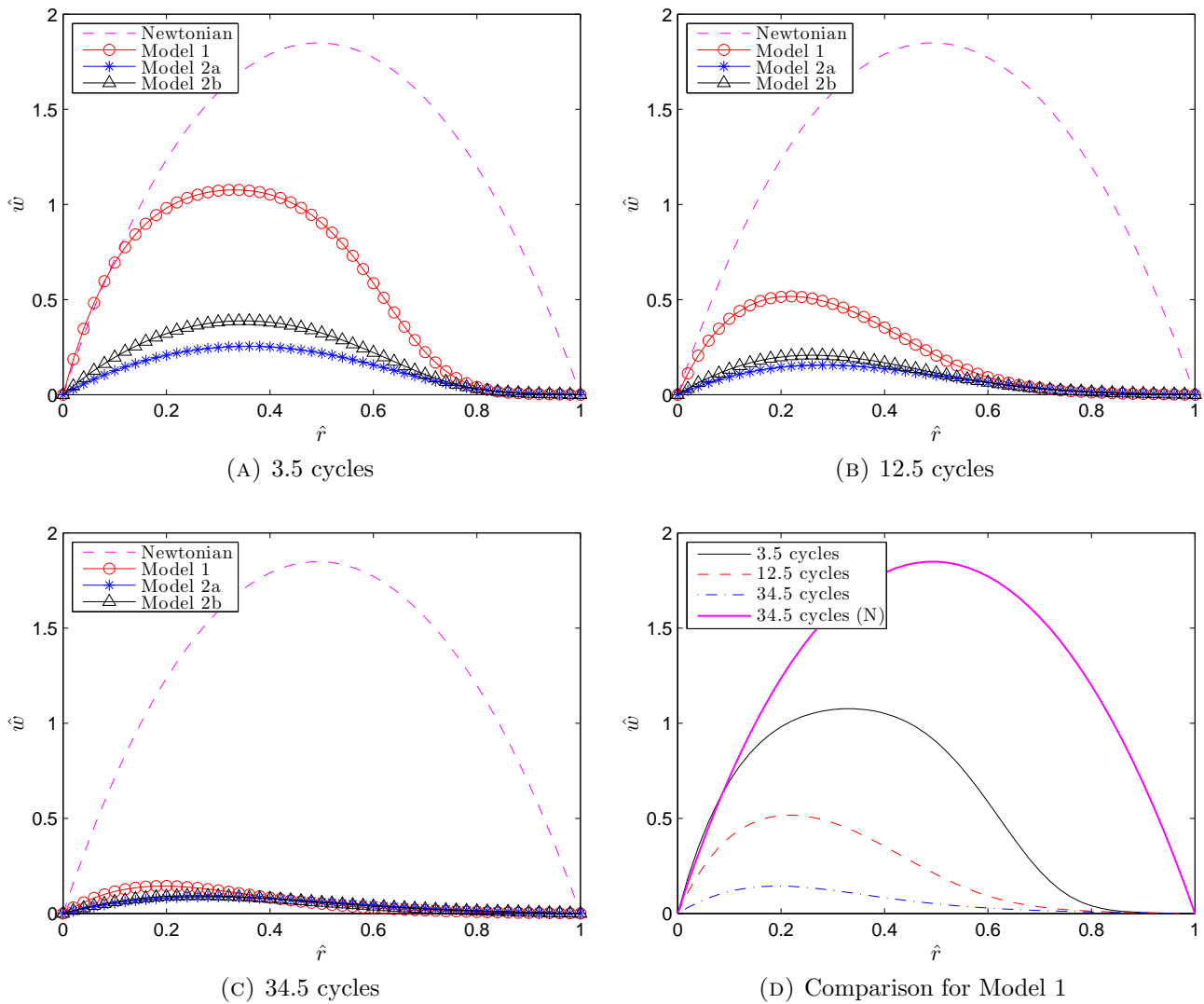

(A) 3.5 cycles

(B) 12.5 cycles

(C) 34.5 cycles

(D) Comparison for Model 1

FIGURE 8. Non-dimensional velocity profiles in the $z$ direction. Comparison for different models at different cycles are shown in (a), (b), (c). In (d), the non-dimensional velocity profiles in the $z$ direction are compared for different cycles for Model 1. The parameters chosen are $r_i = 1$, $r_o = 1.2$, $p_A = -p_B = 1$, $p_f = 1$, $\overline{\omega}_\theta = 1$, $Re = 10$, $Pe = 1000$.

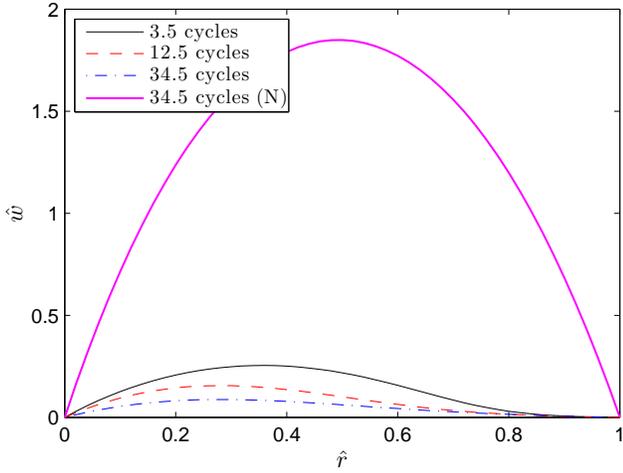

(A) Comparison for Model 2a

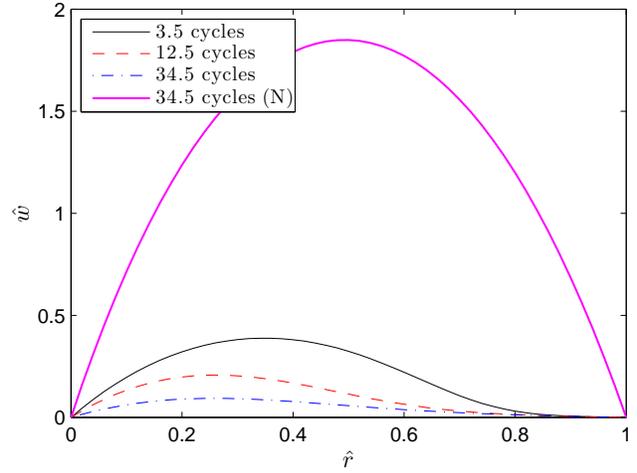

(B) Comparison for Model 2b

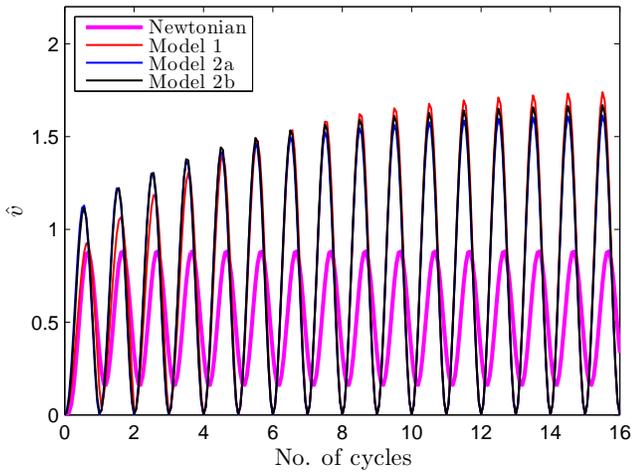

(C) Centerline $\hat{v}$ for different models

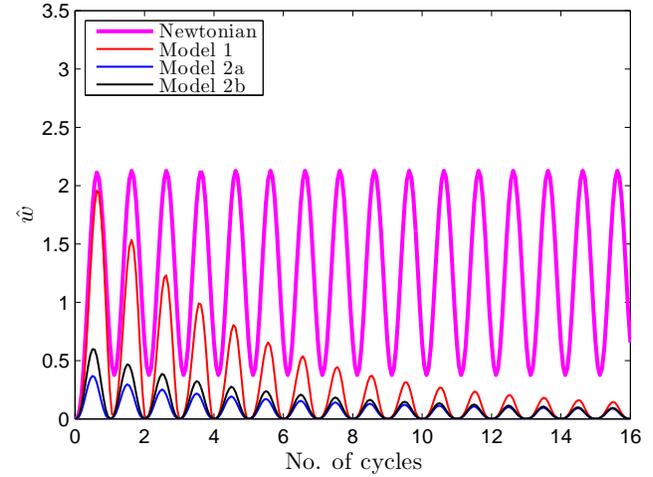

(D) Centerline $\hat{w}$ for different models

FIGURE 9. The non-dimensional velocity profiles in the $z$ direction are compared for different cycles for Model 2a and Model 2b. Results are shown in (a) and (b) respectively. Non-dimensional centerline velocity (at $\hat{r} = 0.5$) in the $\theta$ direction as a function of number of cycles for different models is shown in (c). Non-dimensional centerline velocity in the $z$ direction as a function of number of cycles for different models is shown in (d). The parameters chosen are $r_i = 1$, $r_o = 1.2$, $p_A = -p_B = 1$, $p_f = 1$, $\overline{\omega}_\theta = 1$, $Re = 10$, $Pe = 1000$.

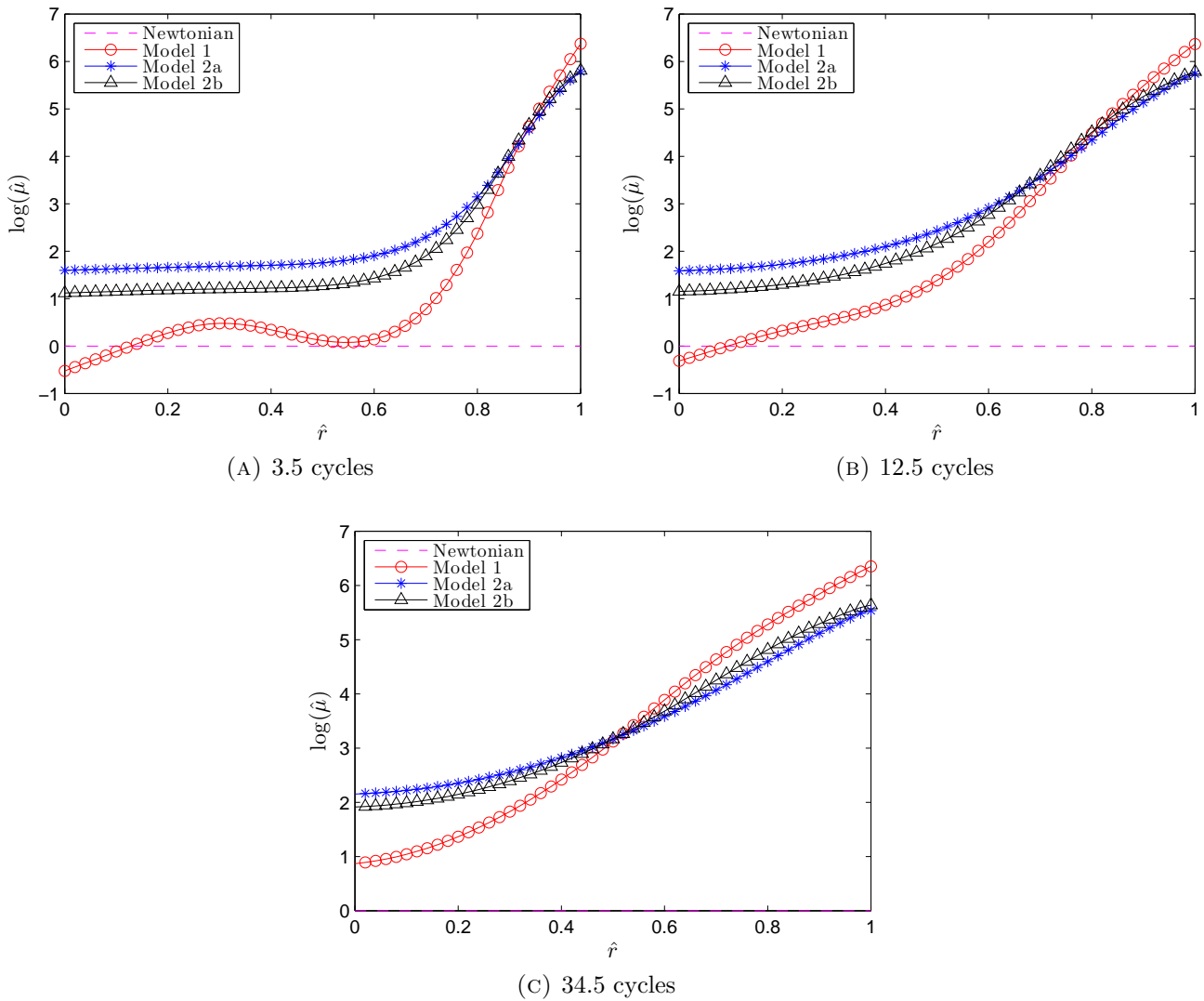

(A) 3.5 cycles

(B) 12.5 cycles

(C) 34.5 cycles

Figure 10. Non-dimensional apparent viscosity profile comparison for different models at different cycles are shown in (a), (b), (c). The parameters chosen are $r_i = 1$, $r_o = 1.2$, $p_A = -p_B = 1$, $p_f = 1$, $\overline{\omega}_\theta = 1$, $Re = 10$, $Pe = 1000$.

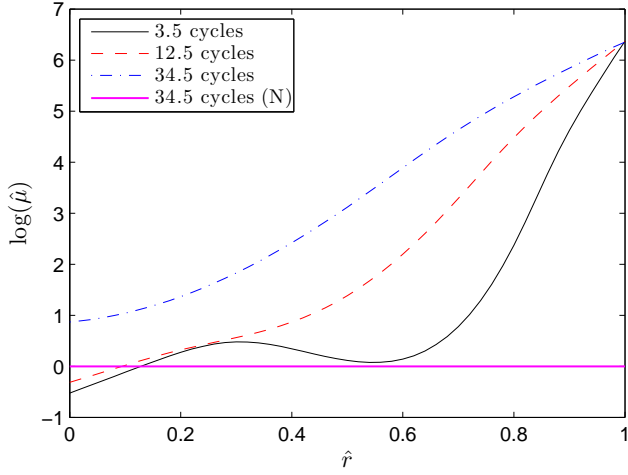

(A) Comparison for Model 1

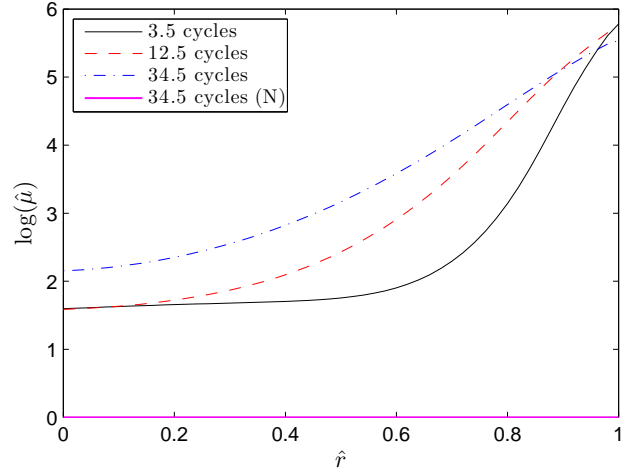

(B) Comparison for Model 2a

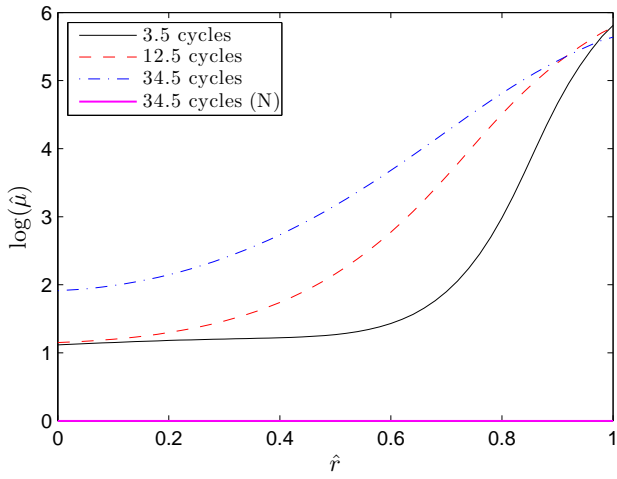

(C) Comparison for Model 2b

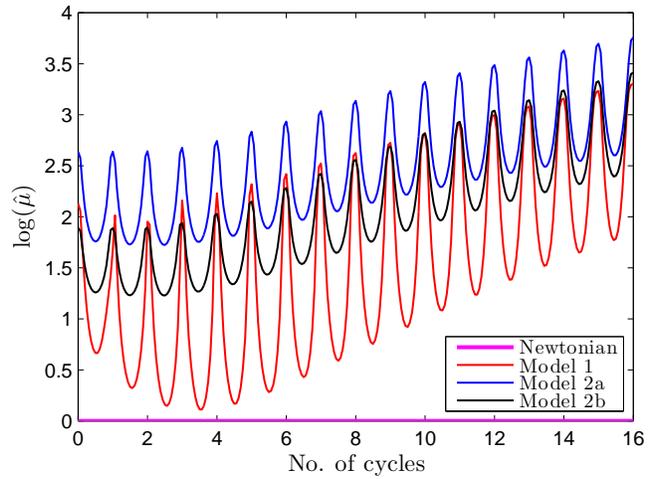

(D) Centerline $\hat{\mu}$ for different models

FIGURE 11. The non-dimensional apparent viscosity profiles are compared for different cycles for Model 1, Model 2a and Model 2b. Results are shown in (a), (b) and (c) respectively. Non-dimensional centerline apparent viscosity (at $\hat{r} = 0.5$) as a function of number of cycles for different models is shown in (d). The parameters chosen are $r_i = 1$, $r_o = 1.2$, $p_A = -p_B = 1$, $p_f = 1$, $\overline{\omega}_\theta = 1$, $Re = 10$, $Pe = 1000$.

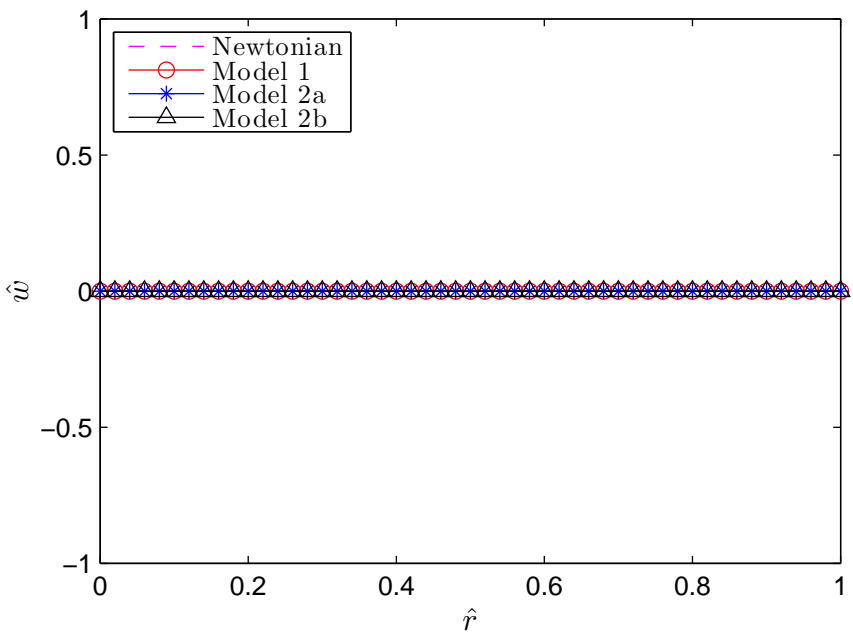